\documentclass[11pt]{article}

\usepackage{epsfig}
\usepackage{graphicx}
\usepackage{latexsym}
\usepackage{amsmath}
\usepackage{amsbsy}
\usepackage{amssymb}
\usepackage{theorem}
\usepackage{amsfonts}

\usepackage{mathenv}

\usepackage[francais,english]{babel}
\usepackage[latin1]{inputenc}
\usepackage[T1]{fontenc}

\def\be{\begin{equation} \displaystyle}
\def\ee{\end{equation} }

\def\bea{\begin{eqnarray}}
\def\eea{\end{eqnarray} }
\def\bean{\begin{eqnarray*}}
\def\eean{\end{eqnarray*} }

\def\div{{\rm div \;}}

\def\R{{\rm I\hspace{-0.50ex}R} }
\def\C{\rm \hbox{C\kern-.57em\raise.47ex
         \hbox{$\scriptscriptstyle |$}\kern+0.5 em }}
\def\E{{\rm I\hspace{-0.50ex}E} } %esp\'erance
\def\tr{\mathrm{tr}\,} %trace d'une matrice
\def\I{\mbox{Id}} % l'identité
\newtheorem{lem}{Lemma}
\newtheorem{rem}{Remark}

\newcommand{\uu}[1]  {{\boldsymbol #1} }
\newcommand{\dps}{\displaystyle   }

\def\XX{\uu{X}}

\def\WW{\uu{W}}
\def\We{\mbox{We}}
\def\Rey{\mbox{Re}}

\newcommand{\U}      { \uu{u}         }

\title{New entropy estimates for the Oldroyd-B model, and related models}
\author{D. Hu$^{(1)}$, T. Lelièvre$^{(2,3)}$\\
\footnotesize{(1) School of Mathematical Sciences, Peking University, Beijing, 100871,
P.R. China.}\\
\footnotesize{(2) CERMICS, Ecole Nationale des
  Ponts (ParisTech), 6 \& 8 Av. B. Pascal,
  77455 Marne-la-Vall\'ee, France.}\\
\footnotesize{(3) INRIA Rocquencourt, MICMAC team, B.P. 105, 78153 Le Chesnay Cedex, France.}\\
 hdmxy@hotmail.com, lelievre@cermics.enpc.fr
}

\begin{document}

\maketitle

\abstract{This short note presents the derivation of a new {\it a priori} estimate for
  the Oldroyd-B model. Such an estimate may provide useful information when investigating the
  long-time behaviour of macro-macro models, and the stability of
  numerical schemes. We show how this estimate can be used as a
  guideline to derive new estimates for other macroscopic models, like
  the FENE-P model.}

\section{Introduction}

We consider the Oldroyd-B model:
\bea
\Rey \left( \frac{\partial \U}{\partial t}+\U.\nabla \U
  \right)=(1- \varepsilon) \Delta \U  - \nabla p + \div \uu{\tau},\label{eq:u}\\
\div (\U) = 0,\label{eq:div_u}\\
\frac{\partial \uu{\tau}}{\partial t} + \U.\nabla \uu{\tau} = \nabla \U  \uu{\tau}
+\uu{\tau}  (\nabla \U)^T - \frac{1}{\We} \uu{\tau} +
\frac{\varepsilon}{\We}\left( \nabla \U + (\nabla \U)^T \right),\label{eq:tau}
\eea
where the Reynolds number $\Rey>0$, the Weissenberg number $\We>0$ and
$\varepsilon \in (0,1)$ are some non-dimensional numbers.
We suppose that the space variable $\uu{x}$ lives in a bounded domain ${\cal D}$
of~$\R^d$. This system is supplied with initial conditions on the
velocity~$\U$ and on the stress tensor~$\uu{\tau}$. For simplicity, we assume no-slip boundary conditions
on the velocity~$\U$:
\be\label{eq:BC}
\U=0 \text{ on $\partial {\cal D}$}.
\ee
We {\em suppose} that the initial data and the geometry are such that there
exists a unique regular solution to~(\ref{eq:u})--(\ref{eq:tau}) and our
aim is to derive some {\it a priori} estimates on this solution.

Let us introduce the so-called conformation tensor
$\uu{A}=\frac{\We}{\varepsilon} \uu{\tau} + \I$. 
The partial differential equation (PDE) on $\uu{\tau}$ translates into the following PDE on~$\uu{A}$:
\begin{equation}\label{eq:A}
\frac{\partial \uu{A}}{\partial t} + \U.\nabla \uu{A} = \nabla \U  \uu{A}
+\uu{A}  (\nabla \U)^T - \frac{1}{\We} \uu{A} + \frac{1}{\We} \I.
\end{equation}
One can check that if
\begin{equation}\label{eq:hyp_A}
\text{$\uu{A}(t=0)=\frac{\We}{\varepsilon} \uu{\tau}(t=0) + \I$ is a  positive definite symmetric matrix,}
\end{equation}
then this property is propagated forward in time by~\eqref{eq:A} (and, in
particular, $\uu{\tau}$ is symmetric). Assuming uniqueness of solution, this can
be proven for example by using the probabilistic
interpretation of $\uu{A}$ as a covariance matrix, as explained in
Section~\ref{sec:entropy}. We will assume throughout this note that~\eqref{eq:hyp_A} is satisfied. Concerning the importance of
positive-definiteness of $\uu{A}$, we refer for example to~\cite[Section
9.8.10]{keunings-89} and also to the recent work~\cite{fattal-kupferman-04,fattal-kupferman-05}.

In Section~\ref{sec:classical}, we recall how the classical {\it a priori}
estimate for the Oldroyd-B model is derived. Next we show how it can be used to
derive some bounds on the stress tensor, provided the initial condition
satisfies $\det \uu{A}(t=0)>1$. In Section~\ref{sec:entropy}, we
establish a new estimate, which comes from an entropy estimate on the
micro-macro model associated with the Oldroyd-B model
(see~\cite{jourdain-le-bris-lelievre-otto-06}). This estimate provides
bounds on $(\U,\uu{\tau})$
 without any assumption on $\uu{\tau}(t=0)$ (apart from~(\ref{eq:hyp_A})). This new estimate could
be useful to study the longtime behaviour of some macro-macro models, or
to analyze the stability of some numerical schemes. Current research is
directed towards clarifying this.

\section{The classical estimate}\label{sec:classical}

Let us first introduce the kinetic energy:
\begin{equation}
E(t)=\frac{1}{2} \int_{{\cal D}} |\U|^2.
\end{equation}
We easily obtain from~(\ref{eq:u})--(\ref{eq:div_u}):
\begin{equation}\label{eq:E_hom}
\Rey \frac{dE}{dt}= - (1-\varepsilon) \int_{{\cal D}} |\nabla \U|^2 -  
\int_{{\cal D}}  \uu{\tau}: \nabla \U,
\end{equation}
where for two matrices $A$ and $B$, we denote $A:B=A_{i,j} B_{i,j}=\tr(A^T
B)$. On the other hand, taking the trace of the PDE~(\ref{eq:tau}) on $\uu{\tau}$ and
integrating over $\cal D$, we
get:
$$
\frac{d}{d t} \int_{{\cal D}} \tr \uu{\tau} = 2 \int_{{\cal D}} \nabla
\U : \uu{\tau} - \frac{1}{\We} \int_{{\cal D}} \tr \uu{\tau}.
$$

We thus obtain the following estimate:
\begin{equation}\label{eq:estim_oldroyd_classical}
\begin{array}{l}
\dps{\frac{d}{dt}\left( \frac{\Rey}{2} \int_{{\cal D}} |\U|^2 +
    \frac{1}{2} \int_{{\cal D}}
    \tr \uu{\tau} \right)} \\
\dps{+ (1-\varepsilon) \int_{{\cal D}} |\nabla \U|^2  +\frac{1}{2
    \We} \int_{\cal D} \tr \uu{\tau} = 0.}
\end{array}
\end{equation}

\begin{rem}
In terms of $\uu{A}$, the energy estimate~(\ref{eq:estim_oldroyd_classical}) writes:
\begin{equation}\label{eq:estim_oldroyd}
\begin{array}{l}
\dps{\frac{d}{dt}\left( \frac{\Rey}{2} \int_{{\cal D}} |\U|^2 +
    \frac{\varepsilon}{2 \We} \int_{{\cal D}}
    \tr \uu{A} \right)} \\
\dps{+ (1-\varepsilon) \int_{{\cal D}} |\nabla \U|^2  +\frac{\varepsilon}{2
    \We^2} \int_{\cal D} \tr (\uu{A} - \I) = 0.}
\end{array}
\end{equation}
\end{rem}

In Lemma~\ref{lem:trace_tau} below, we prove that $\tr \uu{\tau}$ is
positive if $\det \uu{A}(t=0)>1$. This result combined with the
estimate~(\ref{eq:estim_oldroyd_classical}) thus yields some {\it a
  priori} bounds on $(\U,\tau)$ provided $\det(\uu{A})(t=0)>1$.
In particular, it shows that
$\U$ and $\uu{\tau}$ go exponentially fast to $0$ in the long time limit,
using~(\ref{eq:estim_oldroyd_classical}) and the Poincaré inequality:
$\int_{{\cal D}} |\U|^2 \leq C \int_{{\cal D}} |\nabla \U|^2$.

\begin{lem}\label{lem:trace_tau}
Let us assume that $\det \uu{A}(t=0)>1.$ Then, we have $\forall t \geq 0,\,
\det\uu{A} (t)>1$ and this implies that $\tr \uu{\tau}(t)>0.$
\end{lem}
\begin{proof}
Using~(\ref{eq:A}) and the Jacobi identity (which states
that for any invertible matrix $M$ depending smoothly on a parameter $t$,
$\frac{d}{dt} \ln \det M= \tr \left( M^{-1}\frac{dM}{dt}\right)$), we have:
\begin{equation}\label{eq:jacobi}
\frac{\partial \ln(\det \uu{A}) }{\partial t} + \U.\nabla \ln(\det
\uu{A}) = \frac{1}{\We} \tr \left( \uu{A}^{-1}  - \I\right).
\end{equation}
Since for any symmetric positive matrix $M$ of size
$d \times d$, 
\begin{equation}\label{eq:in_mat}
(\det M)^{1/d} \leq (1/d) \tr M,
\end{equation}
 we obtain
\begin{equation*}
\frac{\partial \ln(\det \uu{A}) }{\partial t} + \U.\nabla \ln(\det
\uu{A}) \geq \frac{d}{\We}  \left( (\det \uu{A})^{-1/d}  - 1\right),
\end{equation*}
which we can rewrite in terms of $y=(\det \uu{A})^{1/d}$:
\begin{equation}\label{eq:ineq_jacobi}
\We \left( \frac{\partial y }{\partial t} + \U.\nabla y \right) \geq  \left( 1 - y \right).
\end{equation}
This shows that $y>1$ if $y(t=0)>1$, and thus that $\det \uu{A}>1$ if $\det \uu{A}(t=0)>1$. 

Indeed, using the characteristic method (by integrating the vector field
$\U(t,\uu{x})$), one can
rewrite~(\ref{eq:ineq_jacobi}) as
$$
\We \frac{D y }{D t} \geq  \left( 1 - y \right).
$$
Now, if $y$ does not remain greater than 1, consider the first time $t_0$
such that $y(t_0)=1$. We have on the one hand
$\frac{D y}{D t}(t_0) <0$ and, on the
other hand $\left(1  - y(t_0) \right)=0$. We
reach a contradiction.

We thus have $ \det \uu{A} >1$ and therefore, using
again~(\ref{eq:in_mat}), $\tr \uu{A} > d$. Since
$\uu{\tau}=\frac{\varepsilon}{\We}(\uu{A} - \I)$, this is equivalent to  $\tr \uu{\tau} > 0$. 
\end{proof}

\begin{rem}
If $\det \uu{A}(t=0)<1$ (which is the case if $\tr \uu{\tau}(t=0)<0$),
Equation~(\ref{eq:ineq_jacobi}) shows that $\det \uu{A}$ grows along the
characteristics as long as $\det \uu{A}<1$.
\end{rem}

\section{Entropy estimate}\label{sec:entropy}

We now consider a micro-macro (or multiscale) formulation of the Oldroyd-B model and some
estimates based on entropy, inspired from~\cite{jourdain-le-bris-lelievre-otto-06}.

\subsection{General derivation of the entropy estimate for micro-macro models}

We consider the following system:
\begin{equation}\label{eq:mic_mac}
\left\lbrace
\begin{array}{l}
\displaystyle{\Rey \left( \frac{\partial \U}{\partial t}(t,\uu{x})+\U(t,\uu{x}).\nabla \U(t,\uu{x})
  \right)=(1- \varepsilon) \Delta \U(t,\uu{x})  - \nabla p(t,\uu{x}) + \div \uu{\tau}(t,\uu{x}),}\\
\displaystyle{\div (\U(t,\uu{x})) = 0,}\\
\displaystyle{\uu{\tau}(t,\uu{x})=\frac{\varepsilon}{\We} \left( \int_{\R^d} (\XX \otimes  \nabla
  \Pi (\XX) ) \psi(t,\uu{x},\XX)\,d\XX - \I \right),}\\
\displaystyle{ \frac{\partial \psi}{\partial
  t}(t,\uu{x},\XX)+\U(t,\uu{x}).\nabla_\uu{x} \psi(t,\uu{x},\XX)}\\
\displaystyle{\qquad =-\div_\XX\left( \left(\nabla_\uu{x}
  \U(t,\uu{x})\XX-\frac{1}{2 \We}\nabla \Pi (\XX)\right) \psi(t,\uu{x},\XX)
\right) +\frac{1}{2 \We} \Delta_\XX \psi(t,\uu{x},\XX).}
\end{array}
\right.
\end{equation}
This system is supplied with initial conditions on the
velocity~$\U$ and on the distribution~$\psi$. We recall that we suppose
 no-slip boundary conditions~(\ref{eq:BC}) on the velocity~$\U$. This
system corresponds to a micro-macro model of polymeric fluids, the
polymer being modelled by two beads linked by a spring with potential
energy $\Pi$. The configurational variable $\XX \in \R^d$ models the
end-to-end vector of the polymer. For more details on the modelling, we refer to~\cite{BCAH-87-2,oettinger-95}.

Notice that we could rewrite the former system as a system coupling a PDE and a stochastic
differential equation (SDE), replacing the last two equations by:
\begin{eqnarray}
\lefteqn{ \uu{\tau}(t,\uu{x})= \frac{\varepsilon}{\We} \Big( \E\left( \XX_t(\uu{x}) \otimes  \nabla
  \Pi (\XX_t(\uu{x})) \right) - \I \Big),} \label{eq:tau_sde} \\
\lefteqn{d \XX_t(\uu{x}) + \U(t,\uu{x}).\nabla_\uu{x} \XX_t(\uu{x}) \, dt} \nonumber\\
&&  = \left(\nabla_\uu{x}
  \U(t,\uu{x})\XX_t(\uu{x}) -\frac{1}{2 \We}\nabla \Pi
  (\XX_t(\uu{x}))\right) \, dt + \frac{1}{\sqrt{ \We }} d\WW_t. \label{eq:sde}
\end{eqnarray}
There, $\E$ denotes the expectation, $\WW_t$ denotes a $d$-dimensional standard Brownian motion
independent from the initial condition $(\XX_0(\uu{x}))_{\uu{x} \in {\cal
    D}}$ which is such that, $\forall \uu{x} \in {\cal
    D}$, the law of $\XX_0(\uu{x})$ is
$\psi(0,\uu{x},\XX)\, d\XX$.

Let us introduce the kinetic energy:
\begin{equation}
E(t)=\frac{1}{2} \int_{{\cal D}} |\U|^2.
\end{equation}
We easily obtain:
\begin{equation}
\Rey \frac{dE}{dt}= - (1-\varepsilon) \int_{{\cal D}} |\nabla \U|^2 -  \frac{\varepsilon}{\We}
\int_{{\cal D}} \int_{ \R^d} (\XX\otimes \nabla \Pi(\XX) ) : \nabla
\U \, \psi.
\end{equation}
We now introduce the entropy of the system, namely:
\begin{eqnarray}
H(t)&=& \int_{{\cal D}} \int_{\R^d} \psi(t,\uu{x},\XX) \ln\left( \frac{\psi(t,\uu{x},\XX)}{\psi_\infty(\XX)}\right),  \\
&=&\int_{{\cal D}}\int_{
  \R^d} \Pi(\XX) \psi(t,\uu{x},\XX) + \int_{{\cal D}}\int_{ \R^d} \psi(t,\uu{x},\XX) \ln(\psi(t,\uu{x},\XX)) + C, \nonumber
\end{eqnarray}
with
\begin{equation}
\psi_\infty(\XX) = \frac{ \exp(- \Pi(\XX))} {\int_{\R^d} \exp(- \Pi(\XX))},
\end{equation}
and $C= \ln(\int_{\R^d} \exp(- \Pi(\XX))) |{\cal D}|$. The function $H$
is actually the relative entropy
of $\psi$ with respect to the equilibrium distribution~$\psi_\infty$.

After some computations (see~\cite{jourdain-le-bris-lelievre-otto-06}), we obtain:
\begin{equation}\label{eq:entropy_homogene}
\frac{dH}{dt}= -\frac{1}{2 \We} \int_{{\cal D}}\int_{\R^d} \psi \left| \nabla
 \ln \left( \frac{\psi}{\psi_\infty} \right) \right|^2   + \int_{{\cal D}}\int_{ \R^d} (\XX\otimes \nabla \Pi(\XX) ) : \nabla
\U \, \psi.
\end{equation}

Therefore, introducing the free energy
$F(t)=E(t)+ \frac{\varepsilon}{\We} H(t)$ of the system, we have:
\begin{equation}\label{eq:energie_libre_homogene}
\boxed{
\begin{array}{l}
\dps{\frac{d}{dt}\left( \frac{\Rey}{2} \int_{{\cal D}} |\U|^2 +
    \frac{\varepsilon}{\We} \int_{{\cal D}}
  \int_{\R^d} \psi \ln\left(
    \frac{\psi}{\psi_\infty}\right) \right)} \\
\dps{+ (1-\varepsilon) \int_{{\cal D}} |\nabla \U|^2  +\frac{\varepsilon}{2 \We^2} \int_{{\cal D}} \int_{ \R^d} \psi \left| \nabla
 \ln \left( \frac{\psi}{\psi_\infty} \right) \right|^2=0.}
\end{array}
}
\end{equation}

Using a logarithmic Sobolev inequality with respect to $\psi_\infty$ and
a Poincaré inequality for $\U \in H^1_0({\mathcal D})$, one can then
obtain exponential convergence to equilibrium $\lim_{t \to \infty}(\U,\psi)=(0,\psi_\infty)$
(see~\cite{jourdain-le-bris-lelievre-otto-06}). For some generalizations
to the case $\U \neq 0$ on $\partial {\mathcal D}$, we refer to~\cite{jourdain-le-bris-lelievre-otto-06}.

\subsection{The Oldroyd-B case}

Let us consider the Hookean dumbbell model, for which the potential
$\Pi$ of the entropic force is:
\be\label{eq:hook}
\dps{\Pi(\XX)=\frac{||\XX||^2}{2}}.
\ee

By It\^o's calculus, it is easy to derive from~(\ref{eq:sde})  that $\uu{A}=\E(\XX_t
\otimes \XX_t)$ satisfies the following PDE:
\begin{equation}\label{eq:oldroyd}
\frac{\partial \uu{A}}{\partial t} + \U.\nabla \uu{A} = \nabla \U  \uu{A}
+\uu{A}  (\nabla \U)^T - \frac{1}{\We} \uu{A} + \frac{1}{\We} \I.
\end{equation}
This translates into the following PDE for $\uu{\tau}=
\frac{\varepsilon}{\We}(\uu{A}- \I)$:
\begin{equation}\label{eq:oldroyd_tau}
\frac{\partial \uu{\tau}}{\partial t} + \U.\nabla \uu{\tau} = \nabla \U  \uu{\tau}
+\uu{\tau}  (\nabla \U)^T - \frac{1}{\We} \uu{\tau} +
\frac{\varepsilon}{\We}\left( \nabla \U + (\nabla \U)^T \right).
\end{equation}
The Hookean dumbbell model is thus equivalent to the
Oldroyd-B model (at least for regular enough solutions).

If $\psi(0,\uu{x},.)$ is Gaussian (with zero mean), so is
$\psi(t,\uu{x},.)$:
$$\psi(t,\uu{x},\XX)=\frac{1}{(2\pi)^{d/2} \sqrt{\det(\uu{A})}}
\exp\left(-\frac{\XX^T \uu{A}^{-1}\XX}{2} \right)$$
where $\uu{A}=\E(\XX_t \otimes \XX_t)=\int_{\R^d} \XX \otimes \XX \,
\psi(t,\uu{x},\XX) \, d\XX$ denotes as above the covariance matrix of~$\XX_t$, which
depends on time and also on the space variable $\uu{x}$. The covariance matrix $\uu{A}$ is symmetric and nonnegative. Moreover, since for
  almost all $t \geq 0$, $\int_{\cal D} \int_{\R^d} \psi(t,\uu{x},\XX) \ln \left(
  \frac{\psi(t,\uu{x},\XX)}{\psi_\infty(\XX)} \right) < \infty$, then for
  almost all $t \geq 0$ and for almost all $\uu{x} \in {\cal D}$, $\uu{A}$ is
  positive.

The following explicit expression of the relative entropy can
then be derived:
$$\int_{\cal D}\int_{\R^d} \psi(t,\uu{x},\XX) \ln \left(
  \frac{\psi(t,\uu{x},\XX)}{\psi_\infty(\XX)} \right)\, d
  \XX=\int_{\cal D} \frac{1}{2}\left(- \ln(\det \uu{A}) - d + \tr \uu{A} \right).$$

On the other hand,
$$
\int_{\cal D} \int_{\R^d} \psi(t,\uu{x},\XX) \left| \nabla_{\XX} \ln \left(
  \frac{\psi(t,\uu{x},\XX)}{\psi_\infty(\XX)}\right) \right|^2 \, d
  \XX = \int_{\cal D} \tr ((\I-\uu{A}^{-1})^2 \uu{A}).
$$

Rewriting~(\ref{eq:energie_libre_homogene}), we thus obtain the
following estimate, in terms  of $\uu{A}$:
\begin{equation}\label{eq:energie_libre_homogene_oldroyd}
\boxed{\begin{array}{l}
\dps{\frac{d}{dt}\left( \frac{\Rey}{2} \int_{{\cal D}} |\U|^2 +
    \frac{\varepsilon}{2\We} \int_{{\cal D}}
   \left(- \ln(\det \uu{A} ) - d + \tr \uu{A} \right) \right)} \\
\dps{+ (1-\varepsilon) \int_{{\cal D}} |\nabla \U|^2  +\frac{\varepsilon}{2
    \We^2} \int_{\cal D} \tr ((\I-\uu{A}^{-1})^2\uu{A} ) = 0.}
\end{array}
}
\end{equation}
This is, in the specific case of Hookean dumbbells (that is Oldroyd-B
model) the macroscopic version of~\eqref{eq:energie_libre_homogene}.

Since $- \ln(\det(\uu{A})) - d + \tr(\uu{A}) \geq 0$, this energy
estimate yields some {\it a priori} bounds on $(\U,\uu{A})$, and thus on
$(\U,\uu{\tau})$. In sharp contrast to
the classical estimate~\eqref{eq:estim_oldroyd_classical}, it provides
bounds on $(\U,\uu{\tau})$ without any assumption on
$\uu{\tau}(t=0)$ (apart from~(\ref{eq:hyp_A})). Using a Poincaré inequality and the fact\footnote{which can be seen as the logarithmic Sobolev inequality for
  Gaussian random variables translated on their covariance matrices} that,  for any symmetric positive matrix $M$ of size
$d \times d$, 
$$- \ln(\det M) - d + \tr M \leq \tr ((\I-M^{-1})^2 M )$$
 exponential
convergence to equilibrium ($\lim_{t \to \infty}(\U,\uu{A})=(0,\I)$) can
be obtained from~(\ref{eq:energie_libre_homogene_oldroyd}).

\begin{rem}
Notice that (\ref{eq:energie_libre_homogene_oldroyd}) can be
schematically obtained
as~(\ref{eq:estim_oldroyd})$\displaystyle{-\frac{\varepsilon}{2 \We} \int_{\mathcal D}}$~(\ref{eq:jacobi}).
\end{rem}

\begin{rem}
If $\psi(0,\uu{x},.)$ is not Gaussian, it is always possible to replace
it by a Gaussian initial condition with the same mean and variance, so that
the macroscopic quantities $(\U,p,\uu{A})$ would be the same for the two
initial conditions.
\end{rem}

\subsection{Application to related macroscopic models}

The energy estimate~\eqref{eq:energie_libre_homogene_oldroyd} can be
used as a guideline to derive energy estimates for other macroscopic
models, even though they cannot be recast as a microscopic model of the
form~\eqref{eq:mic_mac}.

Let us consider the example of the FENE-P model~\cite{peterlin-66,bird-dotson-johnson-80}, for which
\begin{align}
\uu{\tau}&=\frac{\varepsilon}{\We} \left( \frac{ \uu{A} }{1 -
    \tr(\uu{A})/b} - \I \right),\\
\frac{\partial \uu{A}}{\partial t} + \U.\nabla \uu{A} &= \nabla \U  \uu{A}
+\uu{A}  (\nabla \U)^T - \frac{1}{\We} \frac{ \uu{A} }{1 -
    \tr(\uu{A})/b} + \frac{1}{\We} \I.\label{eq:FENE-P}
\end{align}
For this model, we assume~\eqref{eq:hyp_A}, and also that
$\tr(\uu{A})(t=0)< b$, and this property is propagated forward in time by~\eqref{eq:FENE-P} (see~\cite{jourdain-lelievre-04}).

Using the same ideas as for the Oldroyd-B model, we consider the
``entropy'' $H(t)=- \ln(\det \uu{A}) - b \ln \left( 1 -
  \tr(\uu{A})/b\right)$, and we compute its time-derivative:
\begin{align*}
\frac{d}{dt} \int_{\mathcal D} - b \ln \left( 1 - \tr(\uu{A})/b\right) &=
2  \int_{\mathcal D} \frac{\nabla \U : \uu{A}}{1 - \tr(\uu{A})/b} +
\frac{1}{\We} \int_{\mathcal D} \left(- \frac{ \tr (\uu{A}) }{(1 -
    \tr(\uu{A})/b)^2} + \frac{d}{1 -
    \tr(\uu{A})/b} \right),\\
\frac{d}{dt} \int_{\mathcal D} \ln ( \det(\uu{A})) &=\frac{1}{\We}
\int_{\mathcal D} \left(- \frac{d}{1 -  \tr(\uu{A})/b} +
  \tr(\uu{A}^{-1}) \right).
\end{align*}
Combining these expressions with~\eqref{eq:E_hom}, we obtain
\begin{equation}\label{eq:E_FENE-P}
\boxed{
\begin{array}{l}
\dps{\frac{d}{dt}\left( \frac{\Rey}{2} \int_{{\cal D}} |\U|^2 +
    \frac{\varepsilon}{2\We} \int_{{\cal D}}
   \left(- \ln(\det \uu{A}) - b \ln \left( 1 - \tr(\uu{A})/b\right) \right) \right)} \\
\dps{+ (1-\varepsilon) \int_{{\cal D}} |\nabla \U|^2  +\frac{\varepsilon}{2
    \We^2} \int_{\cal D}\left(\frac{ \tr (\uu{A}) }{(1 -
    \tr(\uu{A})/b)^2} - \frac{2d}{1 -  \tr(\uu{A})/b} +
  \tr(\uu{A}^{-1})\right)   = 0.}
\end{array}
}
\end{equation}
One can check that for any symmetric positive matrix $M$ of size
$d \times d$:
\begin{equation}\label{eq:H}
- \ln(\det(M)) - b \ln \left( 1 - \tr(M)/b\right) \geq - (b+d) \ln
\left( \frac{b}{b+d} \right)  \geq d 
\end{equation}
and that
\begin{align}
- \ln(\det(M)) &- b \ln \left( 1 - \tr(M)/b\right) + (b+d) \ln
\left( \frac{b}{b+d} \right) \nonumber \\
&\leq \left(\frac{ \tr (M) }{(1 -
    \tr(M)/b)^2} - \frac{2d}{1 -  \tr(M)/b} +
  \tr(M^{-1})\right). \label{eq:FENE-P_HI}
\end{align}
The proof of these inequalities is tedious and can be done by
diagonalizing the matrix~$M$.

Equation~(\ref{eq:H}) shows that $$\frac{\Rey}{2} \int_{{\cal D}} |\U|^2 +
    \frac{\varepsilon}{2\We} \int_{{\cal D}}
   \left(- \ln(\det \uu{A}) - b \ln \left( 1 - \tr(\uu{A})/b\right) + (b+d) \ln
\left( \frac{b}{b+d} \right) \right)$$
is a non-negative quantity, and thus that~\eqref{eq:E_FENE-P} indeed
yields some {\it a priori} bounds on~$(\U,\uu{A})$.

Equation~(\ref{eq:FENE-P_HI}) (which plays the role of
the log-Sobolev inequality in the micro-macro models) shows that the
estimate~\eqref{eq:E_FENE-P} can be used to prove exponential
convergence to equilibrium.

\bibliography{biblio_HD,ma_biblio}
\bibliographystyle{plain}

\end{document}